\documentclass[draft]{amsart}

\usepackage{amssymb}
\usepackage{url}
\usepackage{amscd}

\newtheorem{theorem}{T{\hskip 0pt\footnotesize\bf HEOREM}}[section]
\newtheorem{lemma}[theorem]{L{\hskip 0pt\footnotesize\bf EMMA}}
\newtheorem{proposition}[theorem]{P{\hskip 0pt\footnotesize\bf ROPOSITION}}
\newtheorem{definition}[theorem]{D{\hskip 0pt\footnotesize\bf EFINITION}}
\newtheorem{corollary}[theorem]{C{\hskip 0pt\footnotesize\bf OROLLARY}}

\newtheorem{example}[theorem]{E{\hskip 0pt\footnotesize\bf XAMPLE}}

\newcommand{\Proof}{\noindent{\bf P{\footnotesize \small\bf roof}:}}
\newcommand{\ProofEnd}{\hfill $\Box$ \bigskip}

\def\a {\alpha }

\newcommand {\C}  {{\mathbb C}}

\renewcommand {\phi} {{\varphi}}

\newcommand {\Scirc} {\raise.2ex\hbox{$\scriptstyle\circ$}}

\newcommand{\Ba}[1]{\begin{array}{#1}}
\newcommand{\Ea}{\end{array}}
\newcommand{\Be}{\begin{equation}}
\newcommand{\Ee}{\end{equation}}
\newcommand{\Bea}{\begin{eqnarray}}
\newcommand{\Eea}{\end{eqnarray}}
\newcommand{\Beas}{\begin{eqnarray*}}
\newcommand{\Eeas}{\end{eqnarray*}}
\newcommand{\Benu}{\begin{enumerate}}
\newcommand{\Eenu}{\end{enumerate}}
\newcommand{\Bi}{\begin{itemize}}
\newcommand{\Ei}{\end{itemize}}

\newcommand{\BR}{\begin{Remark} \em}
\newcommand{\ER}{\end{Remark}}
\newcommand{\BE}{\begin{example} \em}
\newcommand{\EE}{\end{example}}

scaled \magstep1

\newcommand {\ub} {{\mathbb B^n}}

\newcommand{{\tBox}}{{\widetilde{\raisebox{-0.2ex}[1.25ex][0ex]{$\Box$}}}}

\newcommand{\bprop} {\begin{proposition}}
\newcommand{\eprop} {\end{proposition}}
\newcommand{\btheo} {\begin{theorem}}
\newcommand{\etheo} {\end{theorem}}
\newcommand{\blem} {\begin{lemma}}
\newcommand{\elem} {\end{lemma}}
\newcommand{\bcor} {\begin{corollary}}
\newcommand{\ecor} {\end{corollary}}

\def\a{\alpha}

\def\B{\mathcal B}

\scrollmode

\begin{document}


\title[$\rho$- Carleson measures on the unit ball]{On some equivalent definitions of $\rho$- Carleson measures on the unit ball}

\author{Benoit F. Sehba}
\address{Beno\^it F. Sehba\\ Department of Mathematics, University of Ghana,\\ P. O. Box LG 62 Legon, Accra, Ghana}
\email{bfsehba@ug.edu.gh}





\maketitle \noindent  
\begin{abstract}
We give in this paper
some equivalent definitions of the so called
$\rho$-Carleson measures when
$\rho(t)=(\log(4/t))^p(\log\log(e^4/t))^q$, $0\le
p,q<\infty$. As applications, we characterize the pointwise
multipliers on $LMOA(\mathbb S^n)$ and from this space to
$BMOA(\mathbb S^n)$. Boundedness of the Ces\`aro type
integral operators on $LMOA(\mathbb S^n)$ and from
$LMOA(\mathbb S^n)$ to $BMOA(\mathbb S^n)$ is considered as
well. This is Chapter 2 of \cite{S2} and some of results were already published
in \cite{BGS} and \cite{S1}.
\end{abstract}

\section{ Holomorphic function spaces and Carleson measures in the unit ball}

\subsection{Some holomorphic function spaces in the unit ball of $\C^n$}
We define here various holomorphic function spaces involved
in this paper. We refer to the book \cite{KZ} for the
proof of different assertions stated below.




Recall that for $\alpha >-1$ the weighted Lebesgue measure
$dV_\alpha$ is defined by \Be dV_\a(z)=c_\a(1-|z|^2)^\a
dV(z),\Ee where \Be
c_\a=\frac{\Gamma(n+\a+1)}{n!\Gamma(\a+1)}\Ee is a
normalizing constant so that $V_\a(\ub)=1$.
\begin{definition} For $\a > -1$ and $0 < p < \infty$,
the weighted Bergman space\index{weighted Bergman space}
$A_{\a}^{p}(\mathbb B^n)$ consists of holomorphic functions
$f$ in $L^p(\mathbb B^n, dV_\a)$, that is \Be
A_{\a}^{p}(\mathbb B^n)=L^p(\mathbb B^n, dV_\a)\bigcap
H(\mathbb B^n).\Ee \end{definition}

We use the notation \Be ||f||_{p,\a}^p:= \int_{\mathbb
B^n}|f(z)|^{p}dV_{\a}(z) \Ee for $f\in L^p(\mathbb B^n,
dV_\a)$.
 \begin{definition}For $0 < p < \infty$ the Hardy space $\mathcal H^p(\mathbb B^n)$ \index{Hardy space} is
 the space of all $f\in H(\mathbb B^n)$ such that \Be ||f||_{p}^{p} := \sup_{0<r<1}\int_{\mathbb S^n}|f(r\xi)|^{p}d\sigma(\xi) <
 \infty.\Ee\end{definition}

 The space of all bounded holomorphic functions in $\mathbb
 B^n$ will be denoted $\mathcal H^\infty(\mathbb B^n)$.

 For any $\xi\in
 \mathbb S^n$ and $\delta > 0$, let $$ B_{\delta}(\xi)=\{w\in \mathbb S^n : |1-\langle w,\xi \rangle |<
 \delta\},$$ and $$Q_{\delta}(\xi)=\{z\in \mathbb B^n : |1-\langle z,\xi \rangle |<
 \delta\}.$$ These are the higher dimension analogues of Carleson
 regions. For $f\in \mathcal H^1(\mathbb B^n)$, we still denote $f(\xi)$, for $\xi\in \mathbb S^n$,
 the admissible limit at the boundary, which exists a.e. \begin{definition}The space of functions of bounded mean
 oscillation\index{BMOA}
 $BMOA$ is the space of all $f\in \mathcal H^1(\mathbb B^n)$ for which there exists a constant $C>0$ so that
$$\sup_{B=B_\delta(\xi), \atop \delta\in ]0,1[,\xi\in \mathbb S^n} \frac{1}{\sigma(B)}\int _B|f-f_{B}|d\sigma\le C.$$ Here and anywhere else,
$f_B$ denotes the mean-value of $f$ on $B$.\end{definition}

The space $BMOA$ is Banach space when equipped with the
norm
$$||f||_{BMOA}^2=|f(0)|+\sup_{B=B_\delta(\xi), \atop \delta\in ]0,1[,\xi\in \mathbb S^n} \frac{1}{\sigma(B)}\int _B|f-f_{B}|d\sigma.$$

We now define the space of functions of logarithmic mean
oscillation\index{LMOA} $LMOA$.
\begin{definition} An analytic function
 $f$ belongs to $LMOA$ if $f\in \mathcal H^1(\mathbb B^n)$
and there exists a constant $C>0$ so that
$$\sup_{B=B_\delta(\xi),\atop \delta\in]0,1[,\xi\in \mathbb S^n} \frac {\log\frac{4}{\delta}}{\sigma(B)}
\int _B|f-f_{B}|d\sigma\le C.$$\end{definition}

The space $LMOA$ is Banach space when equipped with the
norm
$$||f||_{LMOA}^2=|f(0)|+\sup_{B=B_\delta(\xi),\atop \delta\in]0,1[,\xi\in \mathbb S^n} \frac {\log\frac{4}{\delta}}{\sigma(B)}
\int _B|f-f_{B}|d\sigma.$$

\begin{definition} The Bloch space\index{Bloch space}
$\mathcal{B}$ consists of all $f\in H(\mathbb B^n)$ such
 that \Be ||f||_{\mathcal{B}}= |f(0)|+\sup_{z\in B_n} |Rf(z)|(1-|z|^2) <
 \infty.\Ee\end{definition}

 We also recall the following definition of the
logarithmic (weighted) Bloch space $L\mathcal
B$.\begin{definition} An analytic function $f$ belongs to
$L\mathcal B$ if \Be \sup_{z\in B_n}
|Rf(z)|(1-|z|^2)\log\frac{4}{1-|z|^2} <
 \infty.\Ee\end{definition}

 The natural norm on  $L\mathcal B(\mathbb B^n)$ is given by
\Be ||f||_{L\mathcal{B}}= |f(0)|+\sup_{z\in B_n}
|Rf(z)|(1-|z|^2)\log\frac{4}{1-|z|^2} <
 \infty.\Ee Both $\mathcal B$ and $L\mathcal B$ are also
 Banach when equipped with the norms $|| ||_{\mathcal{B}}$
 and $|| ||_{L\mathcal{B}}$ respectively. Moreover, $BMOA$
 continuously embeds in $\mathcal B$ and $LMOA$ embeds
 continuously in $L\mathcal{B}$.

 \subsection{Carleson measures on the unit ball of $\C^n$}

We recall here the definition of Carleson
measures\index{Carleson measure} and their equivalent in
the unit ball. We also introduce Carleson measures with
weight.
\begin{definition} Let $\mu$ denote a positive
Borel measure on $\mathbb B^n$. Then for $0<s<\infty$, the
measure $\mu$ is called a $s$-Carleson measure, if there is
a finite constant $C>0$ such that for any $\xi \in \mathbb
S^n$ and any $0<\delta<1$, \Be \mu(Q_{\delta}(\xi))\le
C(\sigma(B_\delta(\xi)))^{s}.\Ee\end{definition}

When $s=1,$ $\mu$ is just called Carleson measure.  The
infinimum of all these constants $C$ will be denoted
$||\mu||_s$. We will also use $||\mu||$ to denote
$||\mu||_1$. The following theorem is the higher dimension
version of the theorem of L.Carleson \cite{C} and its
reproducing kernel formulation. \btheo\label{B} For a
positive Borel measure $\mu$ on $\mathbb B^n$, and
$0<p<\infty$, the following are equivalent
\begin{itemize}\item[i)] The measure $\mu$ is a Carleson measure
\item[ii)] There is a constant $C_1 >0$ such that, for all
$f\in \mathcal H^p(\mathbb B^n)$, $$\int_{\mathbb
B^n}|f(z)|^{p}d\mu(z) \le C_{1}||f||_{p}^{p}.$$ \item[iii)]
There is a constant $C_2 >0$ such that, for all $a\in
\mathbb B^n$,
$$\int_{\mathbb B^n}\frac{(1-|a|^2)^{n}}{|1-\langle a,w\rangle |^{2n}}d\mu(w)<C_2.$$\end{itemize}\etheo
We say two positive constant $K_1$ and $K_2$ are
comparable, denoted by $K_1\approx K_2$, if there is a
absolute positive constant $M$ such that $$M^{-1}\le
\frac{K_1}{K_2}\le M.$$ We note that the constants
$C_1,\,C_2$ in Theorem \ref{B} are both comparable to
$||\mu||$. The proof of this theorem can be found in
\cite{KZ}. We also have the following theorem in \cite{KZ}
and \cite{ZZ}. \btheo\label{C} For a positive Borel measure
$\mu$ on $\mathbb B^n$, $s>1$ and $0<p<\infty$, the
following are equivalent\begin{itemize}\item[i)] The
measure $\mu$ is a s-Carleson measure\item[ii)] There is a
constant $K_1
>0$ such that, for all $f\in A_{ns-(n+1)}^P$,
$$\int_{\mathbb B^n}|f(z)|^{p}d\mu(z) \le
K_{1}||f||_{p,ns-(n+1)}^{p}.$$\item[iii)] There is a
constant $K_2
>0$ such that, for all $a\in \mathbb B^n$,
$$\int_{\mathbb B^n}\frac{(1-|a|^2)^{ns}}{|1-\langle a,w\rangle |^{2ns}}d\mu(w)<K_2.$$\end{itemize}\etheo
Here both $K_1$ and $K_2$ are comparable to $||\mu||_s$.

We consider here generalized Carleson type measures with
additional logarithmic terms.

\begin{definition}Let $\mu$ be a positive Borel measure on $\mathbb B^n$
and $0< s<\infty$. For $\rho$ a positive function defined
on $(0,1)$, we say $\mu$ is a $(\rho,s)$- Carleson measure
if there is a constant $C>0$ such that for any $\xi \in
\mathbb S^n$ and $0<\delta <1$,\Be \mu(Q_{\delta}(\xi))\le
C\frac{(\sigma(B_\delta(\xi)))^s}{\rho(\delta)}.\Ee
\end{definition}

When $s=1$, $\mu$ is called a $\rho$-Carleson measure. We
are interested in the particular case
$\rho(t)=\rho_{p,q}(t)=(\log(4/t))^p(\log\log(e^4/t))^q$
with $0\le p,q<\infty$. We remark that the case
$\rho(t)=(\log(4/t))^p$ has been studied in \cite{RZ1}. The
corresponding measures in the latter are called
p-logarithmic s-Carleson measures and when $s=1$ we just
called them p-logarithmic Carleson measures and when $p=2$
and $s=1$ we call them logarithmic Carleson measures, using
the vocabulary of \cite{RZ1}. 



\section{The case of $\rho_{p,q}$- Carleson measures}

\subsection{Some useful results}
We give in this subsection some useful results for the the
proof of Theorem \ref{main1} and Theorem \ref{main2}.
\blem\label{estimates} Let $1<N<\infty$ and
$0<\alpha<\infty$. The following assertions
hold.\begin{itemize} \item[i)] For any $0\le p<\infty$,
there exists a positive constant $C_1$ not depending on $N$
so that
$$I_{N,\alpha,p}=\int_1^N\frac{e^{-\alpha t}dt}{(N-t+2)^p}\le
\frac{C_1}{(N+2)^p}.$$\item[ii)] If $\epsilon_1$ and
$\epsilon_2$ are real so that
$\log(2+\epsilon_1)+\epsilon_2>1$, then for any $0\le
p<\infty$, there exists a positive constant $C_2$ not
depending on $N$ so that
$$J_{N,\alpha,p}=\int_1^N\frac{e^{-\alpha t}dt}{(\log(N-t+2+\epsilon_1)+\epsilon_2)^p}\le
\frac{C_2}{(\log(N+2+\epsilon_1)+\epsilon_2)^p}.$$\end{itemize}\elem\Proof\,
$i)$. A simple change of variables gives the following
equalities
$$I_{N,\alpha,p}=\int_2^{N+1}\frac{e^{-\alpha(N+2-x)}dx}{x^p}=e^{-\alpha(N+2)}\int_2^{N+1}x^{-p}e^{\alpha x}dx.$$
Thus $i)$ can be written as
\begin{equation}\label{ineqequiv1}\int_2^{N+1}x^{-p}e^{\alpha x}dx\le
C_1(N+2)^{-p}e^{\alpha (N+2)}.\end{equation} Let
$f(x)=x^{-p}e^{\alpha x}$. Then $f'(x)=x^{-p-1}e^{\alpha
x}(\alpha x-p)>0$ if $x>\frac{p}{\alpha}$. Since $f(x)$ is
obviously continuous on $[2,\infty)$ and increasing as
$x>\frac{p}{\alpha}$, there is a positive constant $K$ such
that for any $x\in [2,N+1]$, $$f(x)\le
Kf(N+1)=K(N+1)^{-p}e^{\alpha(N+1)}.$$ Integrating by parts
gives \Beas \int_2^{N+1}x^{-p}e^{\alpha x}dx &=&
\left. \frac{1}{\alpha}x^{-p}e^{\alpha x} \right|_2^ {N+1}+\frac{p}{\alpha}\int_2^{N+1}x^{-p-1}e^{\alpha x}dx\\
&\le&
\frac{K}{\alpha}(N+1)^{-p}e^{\alpha(N+1)}+\frac{Kp}{\alpha}(N+1)^{-p-1}e^{\alpha(N+1)}\int_2^{N+1}dx\\
&\le& \frac{K(1+p)}{\alpha}(N+1)^{-p}e^{\alpha(N+1)}\\
&\le& \frac{K'(1+p)}{\alpha}(N+2)^{-p}e^{\alpha(N+2)},\Eeas
where $K'$ is another positive constant, independent of
$N$. Thus (\ref{ineqequiv1}) is true, hence $(i)$ is true.

$ii)$. The proof is similar to the proof of $i)$. Let
$x=N+2+\epsilon_1-t$. Then $t=N+2+\epsilon_1-x$, and
$dt=-dx$. Thus
$$J_{N,\alpha,p}=\int_{2+\epsilon_1}^{N+1+\epsilon_1}\frac{e^{-\alpha(N+2+\epsilon_1-x)}dx}{(\epsilon_2+\log x )^p}=
e^{-\alpha(N+2+\epsilon_1)}\int_{2+\epsilon_1}^{N+1+\epsilon_1}(\epsilon_2+\log
x)^{-p}e^{\alpha x}dx.$$ Thus $ii)$ can be written as
\begin{equation}\label{ineqequiv2}\int_{2+\epsilon_1}^{N+1+\epsilon_1}(\epsilon_2+\log x )^{-p}e^{\alpha x}dx\le
C_2[\epsilon_2+\log(N+2+\epsilon_1)]^{-p}e^{\alpha
(N+2+\epsilon_1)}.\end{equation} Let $g(x)=(\epsilon_2+\log
x)^{-p}e^{\alpha x}$. Then
$$g'(x)=(\epsilon_2+\log x)^{-p-1}e^{\alpha x}\left[\alpha(\epsilon_2+\log x)-\frac{p}{x}\right].$$
Since $$\lim_{x\rightarrow
\infty}\left[\alpha(\epsilon_2+\log
x)-\frac{p}{x}\right]=\infty,$$ we know that there exists a
positive constant $M$ such that $g'(x)>0$ for all $x>M$.
Thus $g$ is continuous on $[2+\epsilon_1,\infty)$ and
increasing whenever $x>M$. Therefore there is a positive
constant $K_1$ such that for any $x\in
[2+\epsilon_1,N+1+\epsilon_1]$, $$g(x)\le
K_1g(N+1+\epsilon_1)=K_1[\epsilon_2+\log
(N+1+\epsilon_1)]^{-p}e^{\alpha (N+1+\epsilon_1)}.$$
Integrating by parts gives \Beas
\int_{2+\epsilon_1}^{N+1+\epsilon_1}(\epsilon_2+\log
x)^{-p}e^{\alpha x}dx &=& \left.
\frac{1}{\alpha}(\epsilon_2+\log x)^{-p}e^{\alpha
x}\right|_{2+\epsilon_1}^{N+1+\epsilon_1}+\frac{p}{\alpha}\int_{2+\epsilon_1}^{N+1+\epsilon_1}(\epsilon_2+\log
x)^{-p-1}e^{\alpha x}x^{-1}dx\\
&\le& \frac{K_1}{\alpha}[\epsilon_2+\log
(N+1+\epsilon_1)]^{-p}e^{\alpha (N+1+\epsilon_1)}\\ & &
+\frac{K_1p}{\alpha}[\epsilon_2+\log
(N+1+\epsilon_1)]^{-p-1}e^{\alpha (N+1+\epsilon_1)}\int_{2+\epsilon_1}^{N+1+\epsilon_1}x^{-1}dx\\
&\le& \frac{K_1(1+p)}{\alpha}[\epsilon_2+\log
(N+1+\epsilon_1)]^{-p}e^{\alpha (N+1+\epsilon_1)}\\
&\le& \frac{K_2(1+p)}{\alpha}[\epsilon_2+\log
(N+2+\epsilon_1)]^{-p}e^{\alpha (N+2+\epsilon_1)},\Eeas
where $K_2$ is another positive constant, independent of
$N$. Thus (\ref{ineqequiv2}) is true, hence $ii)$ is true.
The proof is complete.

\ProofEnd

Let set $$K_a(z)=\frac{(1-|a|^2)^n}{|1-\langle a,z\rangle
|^{2n}}.$$ We have the following two theorems
characterizing $(\rho_{p,q},s)$-Carleson measures in the
unit ball.

\btheo\label{general1} Let $0\le p,q <\infty$ and $0 < s
<\infty.$ Let $\mu$ be a positive Borel measure on $\mathbb
B^n$. Then, $\mu$ is a $(\rho_{p,q},s)$-Carleson measure
with $\rho_{p,q}(t)=(\log(4/t))^p(\log\log(e^4/t))^q$, if
and only if\Be\label{equivprincipale} \sup_{a\in \mathbb
B^n}(\log\frac{4}{1-|a|})^{p}(\log\log\frac{e^4}{1-|a|})^q\int_{\mathbb
B^n}K_{a}^{s}(z)d\mu(z)\le C< \infty.\Ee\etheo\Proof The
ideas of the proof are the same as in the proof of Theorem
2 of \cite{RZ1}. We first suppose that $\mu$ is a
$(\rho_{p,q},s)$-Carleson measure and prove
(\ref{equivprincipale}). For $|a|\le \frac{3}{4}$,
(\ref{equivprincipale}) is obvious since the measure is
necessarily finite. Let $\frac{3}{4} < |a| <1$ and choose
$\xi=a/|a|$. For any nonnegative integer $k$, let $r_k =
2^{k-1}(1-|a|)$, $k=1,2,\cdots,N$ and $N$ is the smallest
integer such that $2^{N-2}(1-|a|)\ge 1$. Thus
\Be\label{ineqlog} \log_{2}\frac{4}{1-|a|}\le N \le
1+\log_{2}\frac{4}{1-|a|}.\Ee Let $E_1=Q_{r_1}$ and $E_k
=Q_{r_k}(\xi)-Q_{r_{k-1}}(\xi)$, $k\ge 2.$ We have
$$\mu(E_k)\le \mu(Q_{r_k}(\xi))\le
\frac{C2^{(k-1)ns}(1-|a|)^{ns}}{(\log\frac{4}{2^{k-1}(1-|a|)})^p(\log\log\frac{e^4}{2^{k-1}(1-|a|)})^q}.$$
Moreover, if $k\ge 2$ and $z\in E_k$, then \Beas |1-\langle
a,z\rangle | &=&
|1-|a|+|a|(1-\langle \xi,z\rangle )|\\ &\ge& -(1-|a|)+|a||1-\langle \xi,z\rangle |\\
&\ge& \frac{3}{4}2^{k-1}(1-|a|)-(1-|a|)\\ &\ge&
2^{k-2}(1-|a|).\Eeas
We also have for $z\in E_1$, $$|1-\langle z,a\rangle |\ge
1-|a|>\frac{1}{2}(1-|a|).$$
Using the above estimates, H\"older's inequality, the
equivalence (\ref{ineqlog}) and Lemma \ref{estimates}, we
obtain \Beas \int_{\mathbb B^n}K_{a}^{s}(z)d\mu(z) &\le&
\frac{C}{(1-|a|)^{ns}}
\sum_{k=1}^{N}\frac{1}{2^{2nks}}\frac{r_{k}^{ns}}{(\log\frac{4}{r_k})^p(\log\log\frac{e^4}{r_k})^q}\\
&\le&
C\sum_{k=1}^{N}\frac{1}{2^{kns}}\frac{1}{(\log\frac{4}{2^{k-1}(1-|a|)})^p(\log\log\frac{e^4}{2^{k-1}(1-|a|)})^q}\\
&\le&
C\sum_{k=1}^{N}\frac{1}{2^{ks}}\frac{1}{(\log\frac{4}{2^{k-1}(1-|a|)})^p(\log\log\frac{e^4}{2^{k-1}(1-|a|)})^q}\\
&\le&
C\int_1^N\frac{1}{2^{ts}}\frac{1}{(\log\frac{4}{2^{t-1}(1-|a|)})^p(\log\log\frac{e^4}{2^{t-1}(1-|a|)})^q}dt\\
&\le&
C\left(\int_1^N\frac{1}{2^{ts}}\frac{1}{(\log\frac{4}{2^{t-1}(1-|a|)})^{p+q}}dt\right)^{\frac{p}{p+q}}
\left(\int_1^N\frac{1}{2^{ts}}\frac{1}{(\log\log\frac{e^4}{2^{t-1}(1-|a|)})^{p+q}}dt\right)^{\frac{q}{p+q}}\\
&\le&
\frac{C}{(\log\frac{4}{1-|a|})^p(\log\log\frac{e^4}{1-|a|})^q}.\Eeas
This prove that (\ref{equivprincipale}) holds.

Now, suppose that (\ref{equivprincipale}) holds. For any
$\xi\in \mathbb S^n$ and $0 <\delta <1$, let
$a=(1-\delta)\xi$. Then $1-|a|=\delta$ and for $z\in
Q_{\delta}(\xi)$, we have $K_{a}(z)\ge
\frac{C}{\sigma(B_\delta(\xi))}$. Thus, we obtain easily
that \Beas\infty &>& C\ge
C(\log\frac{4}{1-|a|})^{p}(\log\log\frac{e^4}{1-|a|})^{q}\int_{\mathbb B^n}K_{a}^{s}(z)d\mu(z)\\
&\ge&
C(\log\frac{4}{\delta})^{p}(\log\log\frac{e^4}{\delta})^{q}\int_{Q_{\delta}(\xi)}K_{a}^{s}(z)d\mu(z)\\
&\ge&
C\frac{(\log\frac{4}{\delta})^{p}(\log\log\frac{e^4}{\delta})^{q}}{(\sigma(B_\delta(\xi)))^s}\mu(Q_{\delta}(\xi)).
\Eeas We conclude that $\mu$ is a $(\rho_{p,q},s)$-Carleson
measure. The proof is complete. \ProofEnd

The following is well-known (see also Lemma \ref{loglog}
below).\blem The following assertions are
satisfied.\begin{itemize}\item[i)]There exists a contant
$C>0$ such that for any $f\in BMOA$,
$$|f(z)|\le C\log(\frac{4}{1-|z|})||f||_{BMOA},\,\,\,z\in
\mathbb B^n.$$\item[ii)] Given $a\in \mathbb B^n$, the
function $f_a(z)=\log(\frac{4}{1-\langle z,a\rangle})$ is
uniformly in $BMOA$.\end{itemize}\elem\blem\label{loglog}
The following assertions are
satisfied.\begin{itemize}\item[i)]There exists a contant
$C>0$ such that for any $f\in LMOA$,
$$|f(z)|\le
C\log\log(\frac{e^4}{1-|z|})||f||_{LMOA},\,\,\,z\in \mathbb
B^n.$$\item[ii)] Given $a\in \mathbb B^n$, the function
$f_a(z)=\log\log(\frac{e^4}{1-\langle z,a\rangle})$ is
uniformly in $LMOA$.\end{itemize}\elem\Proof It is not hard
to see that $LMOA$ is a subspace of $L\mathcal B$ with
$||f||_{L\mathcal B}\le C||f||_{LMOA}$. Thus, we only need
to show that $i)$ holds for any $f\in L\mathcal B$.

For any analytic function $f$ in $\mathbb B^n$, one easily
have that $$f(z)-f(0)=\int_0^1\frac{Rf(tz)}{t}dt$$ for all
$z\in \mathbb B^n$. It follows that there exists a contant
$C>0$ such that for any $f\in L\mathcal B$ and any
$z\in \mathbb B^n$, \Beas |f(z)-f(0)| &=& |\int_0^1\frac{Rf(tz)}{t}dt|\\
&\le&
C||f||_{L\mathcal B}\int_0^1\frac{|z|}{(1-|z|t)\log(\frac{e^4}{1-|z|t})}dt\\
&=& C||f||_{L\mathcal
B}(\log\log(\frac{e^4}{1-|z|})-\log4).\Eeas This prove the
pointwise estimate for all $f\in LMOA$.

Let us now prove that the function
$f_a(z)=\log\log(\frac{e^4}{1-\langle z,a\rangle })$ is
uniformly in $LMOA$ or equivalently, by the
characterization of \cite{WS} that the measure
$d\mu_a(z)=|\nabla f_a(z)|^2(1-|z|^2)dV(z)$ is a
logarithmic-Carleson measure (that is a $\rho$-Carleson
measure with $\rho(t)=\log^2(4/t)$) with uniform bound. For
any $\xi\in \mathbb S^n$ and $0< \delta<1$, we set
$$I=\int_{|1-\langle z,\xi\rangle |<\delta}\frac{1-|z|^2}{|1-\langle a,z\rangle |^2|\log(\frac{e^4}{1-\langle z,a\rangle })|^2}dV(z).$$
We have to show that $I\le
C\frac{\sigma(B_\delta(\xi))}{(\log\frac{4}{\delta})^2}$,
where the constant $C>0$ does not depend on the given $a\in
\mathbb B^n$.

If $|1-\langle a,\xi \rangle |\ge 2\delta$ then, for any
$z\in \mathbb B^n$ such that $|1-\langle z,\xi \rangle
|<\delta$, $|1-\langle a,z\rangle |\ge \delta$. Thus,
$$I\le
\delta^{-2}(\log\frac{e^4}{\delta})^{-2}\int_{|1-\langle
z,\xi \rangle |<\delta}(1-|z|^2)dV(z)\lesssim
\frac{\sigma(B_\delta(\xi))}{(\log\frac{e^4}{\delta})^2}.$$
If $|1-\langle a,\xi \rangle |\le 2\delta$, we obtain
 \Beas I &\lesssim&
 \int_{|1-\langle a,z\rangle |< 3\delta}\frac{(1-|z|^2)}{|1-\langle z,a\rangle |^2\log^2\frac{e^4}{|1-\langle z,a\rangle |}}dV(z)\\
 &\lesssim&  \sum_{j=0}^{\infty}\int_{2^{-j-1}.3\delta\le |1-\langle z,a\rangle|\le 2^{-j}.3\delta}\frac{(1-|z|^2)}{|1-\langle z,a\rangle|^2
 \log^2\frac{e^4}{|1-\langle z,a\rangle|}}dV(z)\\
 &\lesssim& \sum_{j=0}^{\infty}2^{2(j+1)}\delta^{-2}(\log2^j\frac{e^4}{\delta})^{-2}\int_{|1-\langle z,a\rangle|\le 2^{-j}.3\delta}(1-|z|^2)dV(z)\\
 &\lesssim&
 \frac{\delta^n}{(\log\frac{e^4}{\delta})^2}\sum_{j=0}^{\infty}2^{2(j+1)}.2^{-j(n+2)}\lesssim
 \frac{\sigma(B_\delta(\xi))}{(\log\frac{4}{\delta})^2}
.\Eeas The proof is complete.\ProofEnd

\subsection{$\rho_{p,q}$- Carleson measures}

In this subsection, we give and prove several equivalent
definitions of $\rho_{p,q}$- Carleson measures. We first
establish a useful lemma. Let $\varphi_z$ be the involutive
automorphism of $\mathbb B^n$ such that $\varphi_z(0)=z$
and $\varphi_z(z)=0$, we remark that for any $a, b,\,\,
\textrm{and}\,\, z\in \mathbb B^n$,
$$K_{a}(z)\cdot K_{b}(\varphi_{a}(z))=K_{\varphi_{a}(b)}(z)$$ and $$K_{a}(\varphi_{a}(z))\cdot K_{a}(z)=1.$$
\blem\label{changevarCarl}Let $0<s<\infty$ and let $\mu$ be
a positive Borel measure on $\mathbb B^n$. Let
$$d\mu_a(z)=\frac{d\mu(\varphi_{a}(z))}{K_{a}^{s}(z)}.$$  Then
$$\sup_{a\in \mathbb B^n}||\mu_a||_s\approx ||\mu||_s.$$\elem\Proof Using the previous remark, we obtain that
\Beas \int_{\mathbb
B^n}K_{b}^{s}(z)\frac{d\mu(\varphi_{a}(z))}{K_{a}^{s}(z)}
&=& \int_{\mathbb
B^n}K_{b}^{s}(\varphi_{a}(w))\frac{d\mu(w)}{K_{a}^{s}
(\varphi_{a}(w))}\\ &=& \int_{\mathbb
B^n}K_{a}^{s}(w)K_{b}^{s} (\varphi_{a}(w))d\mu(w)\\ &=&
\int_{\mathbb B^n}K_{\varphi_{a}(b)}^{s}(w)d\mu(w).\Eeas
The conclusion follows by taking the supremum over $b\in
\mathbb B^n$ and applying Theorem \ref{general1}.\ProofEnd

Let us now recall the following equivalence for the norm of
elements of BMOA space:
$$||f||_{BMOA}\thickapprox \sup_{a\in \mathbb B^n}||f\circ \varphi_a-f(a)||_p$$
for any $0<p<\infty$ (see \cite{KZ}).\blem\label{general3}
Let $0\le p,q<\infty$ and let $\mu$ be a positive Borel
measure on $\mathbb B^n$. Then the following conditions are
equivalent.
\begin{itemize}\item[i)] There exists a positive constant
$C_1$ such that for any $0<\delta<1$ and any $\xi\in
\mathbb S^n$
$$\mu(Q_\delta(\xi))\le C_1\frac{\sigma(B_\delta(\xi))}{(log\frac{4}{\delta})^p(\log\log\frac{e^4}{\delta})^q}.$$
\item[ii)] There exists a positive constant $C_2$ such that
$$\sup_{a\in \mathbb B^n}(\log\frac{4}{1-|a|})^p(\log\log\frac{e^4}{1-|a|})^q\int_{\mathbb
B^n}K_a(z)d\mu(z)\le C_2<\infty.$$

\item[iii)] There exists a positive constant $C_3$ such
that for any $f\in BMOA$,
$$\sup_{a\in \mathbb B^n}(\log\log\frac{e^4}{1-|a|})^q\int_{\mathbb B^n}|f(z)|^pK_a(z)d\mu(z)\le
C_3||f||^p_{BMOA}.$$\item[iv)] There exists a constant
$C_4>0$ such that for any $f\in BMOA$ and any $g\in LMOA$,
$$\sup_{a\in \mathbb B^n}\int_{\mathbb B^n}|f(z)|^p|g(z)|^qK_a(z)d\mu(z)\le
C_4||f||^p_{BMOA}||g||^q_{LMOA}.$$\end{itemize}\elem \Proof
The equivalence $i)\Leftrightarrow ii)$ follows from
Theorem \ref{general1}. We show that $ii)\Rightarrow
iii)\Rightarrow iv)\Rightarrow i)$.

$ii)\Rightarrow iii)$: We first remark that $ii)$ implies
that $\mu$ is a Carleson measure and so is
$\frac{d\mu(\varphi_a(z))}{K_a(z)}$ for any fixed $a\in
\mathbb B^n$ by Lemma \ref{changevarCarl}.

Now, for any $f\in BMOA$, using H\"older's inequality we
obtain \Beas \int_{\mathbb B^n}|f(z)-f(a)|^pK_a(z)d\mu(z)
&\le& (\int_{\mathbb
B^n}|f(z)-f(a)|^{p+q}K_a(z)d\mu(z))^{\frac{p}{p+q}}(\int_{\mathbb
B^n}K_a(z)d\mu(z))^{\frac{q}{p+q}}\\ &\approx&
(\int_{\mathbb
B^n}|fo\varphi_a(z)-f(a)|^{p+q}\frac{d\mu(\varphi_a(z))}{K_a(z)})^{\frac{p}{p+q}}(\int_{\mathbb
B^n}K_a(z)d\mu(z))^{\frac{q}{p+q}}\\ &\le&
C||\mu||^{p/(p+q)}
||fo\varphi_a-f(a)||_{p+q}^p(\int_{\mathbb
B^n}K_a(z)d\mu(z))^{\frac{q}{p+q}}\\ &\le&
C||\mu||^{p/(p+q)}||f||_{BMOA}^p(\int_{\mathbb
B^n}K_a(z)d\mu(z))^{\frac{q}{p+q}}. \Eeas It follows that
\Beas I_1\le
C||\mu||^{p/(p+q)}||f||_{BMOA}^p\left((\log\log\frac{e^4}{1-|a|})^{p+q}\int_{\mathbb
B^n}K_a(z)d\mu(z)\right)^{\frac{q}{p+q}},\Eeas where
$$I_1=(\log\log\frac{e^4}{1-|a|})^q\int_{\mathbb
B^n}|f(z)-f(a)|^pK_a(z)d\mu(z). $$ It is also clear that
$ii)$ implies that $\mu$ is a $\rho$-Carleson measure with
$$\rho(t)=(\log\log\frac{e^4}{t})^{p+q},\,\,\,\, t\in (0,1),$$ which is
equivalent to saying there exists a constant $C>0$ so that
$$(\log\log\frac{e^4}{1-|a|})^{p+q}\int_{\mathbb
B^n}K_a(z)d\mu(z)\le C<\infty.$$ We conclude that
\Be\label{part1}
I_1=(\log\log\frac{e^4}{1-|a|^2})^q\int_{\mathbb
B^n}|f(z)-f(a)|^pK_a(z)d\mu(z) \le C||f||_{BMOA}^p.\Ee
Since $f\in BMOA$, we already know that there exists $C>0$
so that
$$|f(a)|\le C\log\frac{4}{1-|a|}||f||_{BMOA}.$$ Thus, setting $$I_2=(\log\log\frac{e^4}{1-|a|})^q\int_{\mathbb
B^n}|f(a)|^pK_a(z)d\mu(z),$$ we obtain \Beas I_2\le
C(\log\frac{4}{1-|a|})^p(\log\log\frac{e^4}{1-|a|})^q||f||_{BMOA}^p\int_{\mathbb
B^n}K_a(z)d\mu(z).\Eeas We conclude using Theorem
\ref{general1} that \Be\label{part2}
I_2=(\log\log\frac{e^4}{1-|a|})^q\int_{\mathbb
B^n}|f(a)|^pK_a(z)d\mu(z) \le C||f||_{BMOA}^p,\Ee where $C$
is a constant independent of $a$. Finally, we obtain
combining (\ref{part1}) and (\ref{part2}) that for any
$a\in \mathbb B^n$,
\Beas(\log\log\frac{e^4}{1-|a|})^q\int_{\mathbb
B^n}|f(z)|^pK_a(z)d\mu(z) &\le& 2^p(I_1+I_2)\\
&\le& C_2||f||_{BMOA}^p. \Eeas

$iii)\Rightarrow iv)$: For any $f\in BMOA$, let
$d\mu_f(z)=\frac{|f(z)|^p}{||f||_{BMOA}^p}d\mu(z)$. We
would like to show that $iii)$ implies that there exists a
positive constant $C_4$ such that for any $f\in BMOA$ and
any $g\in LMOA$,
$$\sup_{a\in \mathbb B^n}\int_{\mathbb
B^n}|g(z)|^qK_a(z)d\mu_f(z)\le C_4||g||_{LMOA}^q.$$ We
remark that $iii)$ implies in particular that for any $f\in
BMOA$, the measure $d\mu_f$ is a Carleson measure with
$||\mu_f||\approx ||\mu||$. It follows easily as before
that \Be\label{part3}\int_{\mathbb
B^n}|g(z)-g(a)|^qK_a(z)d\mu_f(z)\le C||\mu||\times
||g||_{BMOA}^q\le C||\mu||\times ||g||_{LMOA}^q.\Ee Now,
using the pointwise estimate for $g\in LMOA$, we obtain
\Beas \int_{\mathbb B^n}|g(a)|^qK_a(z)d\mu_f(z) &\le&
C||g||_{LMOA}^q(\log\log\frac{e^4}{1-|a|})^q\int_{\mathbb
B^n}K_a(z)d\mu_f(z).\Eeas It follows using $iii)$ that
there exists $C>0$ so that \Be\label{part4} \int_{\mathbb
B^n}|g(a)|^qK_a(z)d\mu_f(z) \le C||g||_{LMOA}^q.\Ee
Finally, using inequalities (\ref{part3}) and
(\ref{part4}), we conclude that for any $a\in \mathbb B^n$,
\Beas \int_{\mathbb B^n}|f(z)|^p|g(z)|^qK_a(z)d\mu(z) &\le&
2^q\int_{\mathbb
B^n}|f(z)|^p(|g(z)-g(a)|^q+|g(a)|^q)K_a(z)d\mu(z)\\ &\le&
C_2||f||_{BMOA}^p||g(z)||_{LMOA}^q, \Eeas which is $iv)$.

$iv)\Rightarrow i)$: For any $0<\delta<1$ and $\xi\in
\mathbb S^n$, let $a=(1-\delta)\xi$. From $iv)$, we have in
particular that there exists $C>0$ so that for any $f\in
BMOA$ and any $g\in LMOA$,
$$\int_{Q_\delta(\xi)}|f(z)|^p|g(z)|^qK_a(z)d\mu(z)\le
C||f||_{BMOA}^p||g(z)||_{LMOA}^q.$$ We test the above
inequality with $f(z)=f_a(z)=\log\frac{4}{1-\langle
a,z\rangle }$ and $g(z)=g_a(z)=\log\log\frac{e^4}{1-\langle
a,z\rangle }$ which are uniformly in $BMOA$ and $LMOA$
respectively. Remarking that for $z\in Q_\delta(\xi)$,
$K_a(z)\ge \frac{C}{\sigma(B_\delta(\xi))}$,
$\log\frac{4}{\delta}\le |f_a(z)|$ and
$\log\log\frac{e^4}{\delta}\le |g_a(z)|$, we obtain
\Beas\frac{C}{\sigma(B_\delta(\xi))}(\log\frac{4}{1-|a|})^p(\log\log\frac{e^4}{1-|a|})^q\int_{Q_\delta(\xi)}d\mu(z)
&\le& \int_{Q_\delta(\xi)}|f_a(z)|^p|g_a(z)|^qK_a(z)d\mu(z)\\
&\le& C'<\infty.\Eeas That is $$\mu(Q_\delta(\xi))\le
C\frac{\sigma(B_\delta(\xi))}{(\log\frac{4}{\delta})^p(\log\log\frac{e^4}{\delta})^q}.$$
The proof is complete.\ProofEnd

Taking $q=0$ in the above lemma, we obtain the following
corollary (see also \cite{RZ1}). \bcor\label{Zhaocase} Let
$0\le p<\infty$ and let $\mu$ be a positive Borel measure
on $\mathbb B^n$. Then the following conditions are
equivalent.
\begin{itemize}\item[i)] There exists a positive constant
$C_1$ such that for any $0<\delta<1$ and any $\xi\in
\mathbb S^n$
$$\mu(Q_\delta(\xi))\le C_1\frac{\sigma(B_\delta(\xi))}{(\log\frac{4}{\delta})^p}.$$
\item[ii)] There exists a positive constant $C_2$ such that
for any $f\in BMOA$,
$$\sup_{a\in \mathbb B^n}\int_{\mathbb B^n}|f(z)|^pK_a(z)d\mu(z)\le
C_2||f||^p_{BMOA}.$$\end{itemize}\ecor\blem\label{general4}
Let $0\le p,q<\infty$ and let $\mu$ be a positive Borel
measure on $\mathbb B^n$. Then the following conditions are
equivalent.
\begin{itemize}\item[i)] There exists a positive constant
$C_1$ such that for any $0<\delta<1$ and any $\xi\in
\mathbb S^n$,
$$\mu(Q_\delta(\xi))\le C_1\frac{\sigma(B_\delta(\xi))}{(\log\frac{4}{\delta})^p(\log\log\frac{e^4}{\delta})^q}.$$
\item[ii)] There exists a positive constant $C_2$ such that
for any $g\in LMOA$,
$$\sup_{a\in \mathbb B^n}(\log\frac{4}{1-|a|})^p\int_{\mathbb B^n}|g(z)|^qK_a(z)d\mu(z)\le
C_2||g||^q_{LMOA}.$$\end{itemize}\elem \Proof  By Lemma
\ref{general3}, the assertion $i)$ is equivalent to saying
there exists a constant $C>0$ such that for any $f\in BMOA$
and any $g\in LMOA$,
$$\sup_{a\in \mathbb B^n}\int_{\mathbb B^n}|f(z)|^p|g(z)|^qK_a(z)d\mu(z)\le
C||f||_{BMOA}^p||g(z)||_{LMOA}^q.$$ It follows from
Corollary \ref{Zhaocase} that the latter is equivalent to
saying that there exists a positive constant $C$ such that
$$\sup_{a\in \mathbb B^n}(\log\frac{4}{1-|a|})^p\int_{\mathbb B^n}K_a(z)d\mu_g(z)\le C<\infty,$$
where
$d\mu_g(z)=\frac{|g(z)|^q}{||g(z)||_{LMOA}^q}d\mu(z)$. This
proves $ii)$. The proof is complete.\ProofEnd

\btheo\label{main1} Let $0\le p,q<\infty$ and let $\mu$ be
a positive Borel measure on $\mathbb B^n$. Then the
following conditions are
equivalent.\begin{itemize}\item[i)] There is $C_1
>0$ such that for any $\xi\in \mathbb S^n$ and  $0< \delta <1$, $$\mu(Q_\delta(\xi))\le
C_1\frac{\sigma(B_\delta(\xi))}{(\log\frac{4}{\delta})^p(\log\log\frac{e^4}{\delta})^q}.$$
\item[ii)] There is $C_2 >0$ such that $$\sup_{a\in \mathbb
B^n}(\log\frac{4}{1-|a|})^p(\log\log\frac{e^4}{1-|a|})^q\int_{\mathbb
B^n}K_a(z)d\mu(z)\le C_2 <\infty.$$\item[iii)]There is $C_3
>0$ such that for any $f\in BMOA$,  $$\sup_{a\in \mathbb
B^n}(\log\log\frac{e^4}{1-|a|})^q\int_{B_n}|f(z)|^{p}K_{a}(z)d\mu(z)\le
C_3||f||_{BMOA}^{p}.$$\item[iv)]There is $C_4 >0$ such that
for any $g\in LMOA$,  $$\sup_{a\in \mathbb
B^n}(\log\frac{4}{1-|a|})^p\int_{B_n}|g(z)|^{q}K_{a}(z)d\mu(z)\le
C_4||g||_{LMOA}^{q}.$$\item[v)]There is $C_5 >0$ such that
for any $f\in BMOA$ and any $g\in LMOA$,  $$\sup_{a\in
\mathbb B^n}\int_{\mathbb
B^n}|f(z)|^p|g(z)|^{q}K_{a}(z)d\mu(z)\le
C_5||f||_{BMOA}^{p}||g||_{LMOA}^{q}.$$ 
\item[vi)] For $0<r<\infty$, there is $C_6
>0$ such that for any $f\in BMOA$ and any $g\in LMOA$ and any $h\in \mathcal H^r(\mathbb B^n)$,
$$\int_{\mathbb B^n}|f(z)|^{p}|g(z)|^{q}|h(z)|^rd\mu(z)\le
C_6||f||_{BMOA}^{p}||g||_{LMOA}^{q}||h||_r^r.$$\end{itemize}\etheo
\Proof We already have from Lemma \ref{general3} and Lemma
\ref{general4}  that $i)\Leftrightarrow ii)\Leftrightarrow
iii)\Leftrightarrow iv)\Leftrightarrow v)$.
Let
$$d\mu_{f,g}(z)=\frac{|f(z)|^p|g(z)|^q}{||f(z)||_{BMOA}^p||g(z)||_{LMOA}^q}d\mu(z).$$
Then $v)$ is equivalent to saying that
$$\sup_{a\in \mathbb B^n}\int_{\mathbb B^n}K_a(z)d\mu_{f,g}<C_5.$$ By Theorem \ref{B}, this is equivalent to $vi)$. The proof is complete.
\ProofEnd

\subsection{Some applications of $\rho_{p,q}$- Carleson measures}

As first application of Theorem \ref{main1}, we consider
the Cesaro-type integral operator\index{Cesaro-type
integral operator} $T_b$ defined by
$$T_b(f)(z)=\int_{0}^{1}f(tz)Rb(tz)\frac{dt}{t}, \,\,\,\,\, b,\,\,f\in
H(\mathbb B^n).$$ The characterization of the boundedness
properties of $T_b$ has been considered in \cite{AS1},
\cite{AS2}, \cite{SZ} and \cite{RZ1} for the case of the
unit disc.  We first prove the following
result on the boundedness of $T_b$ on $LMOA$.\bcor For
$b\in H(\mathbb B^n)$, $T_b$ is bounded on $LMOA$ if and
only if \Be\label{T_b1} \sup_{a\in \mathbb
B^n}(\log\frac{4}{1-|a|})^2(\log\log\frac{e^4}{1-|a|})^2\int_{\mathbb
B^n}|Rb(z)|^2(1-|z|^2)K_a(z)dV(z)<\infty.\Ee\ecor\Proof We
know from \cite{WS} that, an analytic $b$ is in $LMOA$ if
and only if $(1-|z|^2)|Rb(z)|^2dV(z)$ is a $\rho$-Carleson
measure with $\rho(t)=(\log(4/t))^2$, which by Lemma
\ref{general1} is equivalent to
$$\sup_{a\in
\mathbb B^n}(\log\frac{4}{1-|a|})^2\int_{\mathbb
B^n}|Rb(z)|^2(1-|z|^2)K_a(z)dV(z)<\infty.$$ It is not hard
to see  that
$$R[T_b(f)](z)=f(z)Rb(z).$$ It follows that $T_b$ is bounded on $LMOA$ if and only if for any $f\in LMOA$,
$$\sup_{a\in\
\mathbb B^n}(\log\frac{4}{1-|a|})^2\int_{\mathbb
B^n}|f(z)|^2|Rb(z)|^2(1-|z|^2)K_a(z)dV(z)<C||f||_{LMOA}^2,$$
which by Theorem \ref{main1} is equivalent to saying that
the measure $|Rb(z)|^2(1-|z|^2)dV(z)$ satisfies
$$\sup_{a\in
\mathbb
B^n}(\log\frac{2}{1-|a|})^2(\log\log\frac{e^4}{1-|a|})^2\int_{\mathbb
B^n}|Rb(z)|^2(1-|z|^2)K_a(z)dV(z)<\infty.$$ The proof is
complete.\ProofEnd\bcor For $b\in H(\mathbb B^n)$, $T_b$ is
bounded from $LMOA$ to $BMOA$ if and only if
\Be\label{T_b2} \sup_{a\in \mathbb
B^n}(\log\log\frac{e^4}{1-|a|})^2\int_{\mathbb
B^n}|Rb(z)|^2(1-|z|^2)K_a(z)dV(z)<\infty.\Ee\ecor\Proof It
is well-known that an analytic $b$ is in $BMOA$ if and only
if $(1-|z|^2)|Rb(z)|^2dV(z)$ is a Carleson measure that is
$$\sup_{a\in
\mathbb B^n}\int_{\mathbb
B^n}|Rb(z)|^2(1-|z|^2)K_a(z)dV(z)<\infty.$$ It follows that
$T_b$ is bounded from $LMOA$ to $BMOA$ if and only if for
any $f\in LMOA$,
$$\sup_{a\in\
\mathbb B^n}\int_{\mathbb
B^n}|f(z)|^2|Rb(z)|^2(1-|z|^2)K_a(z)dV(z)<C||f||_{LMOA}^2$$
which by Theorem \ref{main1} is equivalent to saying that
the measure $|Rb(z)|^2(1-|z|^2)dV(z)$ satisfies
$$\sup_{a\in
\mathbb B^n}(\log\log\frac{e^4}{1-|a|})^2\int_{\mathbb
B^n}|Rb(z)|^2(1-|z|^2)K_a(z)dV(z)<\infty.$$ The proof is
complete.\ProofEnd

We also obtain in the same way the following result.\bcor
For $b\in H(\mathbb B^n)$, $T_b$ is bounded on $BMOA$ if
and only if \Be\label{T_b3} \sup_{a\in \mathbb
B^n}(\log\frac{4}{1-|a|})^2\int_{\mathbb
B^n}|Rb(z)|^2(1-|z|^2)K_a(z)dV(z)<\infty.\Ee\ecor

Our next application is about the pointwise
multipliers\index{pointwise multiplier} on $LMOA$. Given
two Banach spaces of analytic functions $X$ and $Y$, we
denote by $\mathcal{M}(X,Y)$ the space of multipliers from
$X$ to $Y$, that is
$$\mathcal{M}(X,Y)=\{f\in H(\mathbb B^n):\,\,\,f\cdot g\in Y\,\,\textrm{for any}\,\, g\in
X\}.$$ When $X=Y$, we just write $\mathcal M(X,X)=\mathcal
M(X)$. The following lemma is an easy adaptation of
\cite[Lemma 3.20]{KZ}.\blem\label{Mutipcond} Suppose that
$X$ and $Y$ are two Banach spaces of holomorphic functions.
If $X$ contains constant functions and each point
evaluation is a bounded linear functional on $Y$, then
every pointwise multiplier from $X$ to $Y$ is in $\mathcal
H^\infty(\mathbb B^n)$.\elem We have the following
characterization of $\mathcal{M}(LMOA)$ for the unit ball
of $\C^n$\bcor\label{LMOMultip} An analytic function $f$ on
$\mathbb B^n$ belongs to $\mathcal{M}(LMOA)$ if and only if
$f\in \mathcal H^\infty(\mathbb B^n)$ and satisfies
(\ref{T_b1}).\ecor\Proof Instead of using Lemma
\ref{Mutipcond}, we give a direct proof of the fact that
any element in $\mathcal{M}(LMOA)$ is necessarily bounded.
For this, we recall that for any $f\in LMOA$,
$$|f(z)|\le C||f||_{LMOA}\log\log\frac{e^4}{1-|z|^2}.$$ Now, for any $a\in
\mathbb B^n$, let $f_a(z)=\log\log(\frac{e^4}{1-\langle
z,a\rangle })$. $f_{a}\in LMOA$ and $||f_a||_{LMOA}\le C
<\infty$.

It follows from these two facts that, if $f\in
\mathcal{M}(LMOA)$,
 then $f\cdot f_a\in LMOA$ and for any $z\in \mathbb B^n$,
 $$|f(z)f_a(z)|\le C||f\cdot f_a||_{LMOA}\log\log\frac{e^4}{1-|z|^2}.$$
Taking $z=a$ in the above inequality, we obtain
$$|f(a)|\le C||f\cdot f_a||_{LMOA}<\infty.$$ That is $f\in H^\infty(\mathbb B^n)$.

That $f\in \mathcal{M}(LMOA)$ means that for any $g\in
LMOA$, the measure $|R(fg)(z)|^2(1-|z|^2)dV(z)$ is a
logarithmic Carleson measure or equivalently that
\Be\label{equivmultip1}I_f(g)\le C||g||_{LMOA}^2,\Ee where
$$I_f(g)=\sup_{a\in \mathbb
B^n}(\log\frac{4}{1-|a|})^2\int_{\mathbb
B^n}|g(z)Rf(z)+f(z)Rg(z)|^2(1-|z|^2)K_a(z)dV(z).$$ Since
$f\in H^\infty(\mathbb B^n)$ and $|Rg(z)|^2(1-|z|^2)dV(z)$
is a logarithmic Carleson measure,
$$\sup_{a\in \mathbb
B^n}(\log\frac{4}{1-|a|})^2\int_{\mathbb
B^n}|f(z)Rg(z)|^2(1-|z|^2)K_a(z)dV(z)\le
C||f||_{\infty}^2||g||_{LMOA}^2.$$ We deduce that if $f\in
H^\infty(B_n)$, then (\ref{equivmultip1}) is equivalent to
$$\sup_{a\in \mathbb
B^n}(\log\frac{4}{1-|a|})^2\int_{\mathbb
B^n}|g(z)|^2|Rf(z)|^2(1-|z|^2)K_a(z)dV(z)\le
C||g||_{LMOA}^2,$$ which by Theorem \ref{main1} is
equivalent to saying that $|Rf(z)|^2(1-|z|^2)dV(z)$
satisfies $$\sup_{a\in \mathbb
B^n}(\log\frac{2}{1-|a|})^2(\log\log\frac{e^4}{1-|a|})^2\int_{\mathbb
B^n}|Rb(z)|^2(1-|z|^2)K_a(z)dV(z)<\infty.$$ The proof is
complete. \ProofEnd

Similarly, we can prove the following results.
\bcor\label{LMOBMOMultip} An analytic function $f$ on
$\mathbb B^n$ belongs to $\mathcal{M}(LMOA,BMOA)$ if and
only if $f\in \mathcal H^\infty(\mathbb B^n)$ and satisfies
(\ref{T_b2}).\ecor
\bcor\label{BMOMultip} An analytic function $f$ on $\mathbb
B^n$ belongs to\\ $\mathcal{M}(BMOA)$ if and only if $f\in
\mathcal H^\infty(\mathbb B^n)$ and satisfies
(\ref{T_b3}).\ecor

The orthogonal projection of $L^2(\partial \mathbb
 B^n)$ onto $\mathcal H^2(\mathbb B^n)$ is called the
 Szeg\"o projection and denoted $P$. It is given by
 \begin{equation} P(f)(z)=\int_{\partial \mathbb B^n}S(z,\xi)f(\xi)d\sigma(\xi),\end{equation}
 where $S(z,\xi)=\frac{1}{(1-\langle z,\xi \rangle)^n}$ is
 the Szeg\"o kernel on $\partial \mathbb B^n$. We denote as
 well by $P$ its extension to $L^1(\partial \mathbb
 B^n)$.

 For $b\in \mathcal H^2(\mathbb B^n)$, the small Hankel
 operator\index{Hankel operator} with symbol $b$ is defined for $f$ a bounded
 holomorphic function by $h_b(f):=P(b\overline f)$. As last
 application, we prove that if $b\in LMOA$, then the Hankel
 operator $h_b$ is bounded on $\mathcal H^1(\mathbb B^n)$.
\vskip .2cm
The following lemma can be proved using integration by
parts (see \cite{Rud}) .\blem\label{integrationparts} Let
$f,g$ be holomorphic polynomials on $\mathbb B^n$. Then the
following identity holds \Beas \int_{\mathbb
S^n}f(\xi)\overline g(\xi)d\sigma(\xi) &=&
C_1\int_{\mathbb B^n}f(z)\overline
{g(z)}dV(z)+C_2\int_{\mathbb B^n}Rf(z)\overline
{g(z)}(1-|z|^2)dV(z)+\\ & & C_3\int_{\mathbb
B^n}f(z)\overline {Rg(z)}(1-|z|^2)dV(z)\Eeas
$C_1,\,\,\,C_2$ and $C_3$ being constants independent of
$f$ and $g$.\elem 

 \btheo The Hankel operator $h_b$ extends into a
bounded operator on  $\mathcal H^1(\mathbb B^n)$ if $b\in
LMOA$.\etheo\Proof Let $b\in LMOA$ or equivalently, such
that $(1-|z|^2)|\nabla b(z)|^2dV(z)$ is a logarithmic
Carleson measure. For $f\in \mathcal H^1(\mathbb B^n)$ and
$g\in BMOA$, we want to estimate $|<h_b(f),g>| = |<b,fg>|$.
Applying Lemma \ref{integrationparts} to $<b,fg>$, it comes
that we only need to estimate the following three
integrals: \Beas I_1:= \int_{\mathbb B^n}|f(z)||
g(z)||b(z)|dV(z),\Eeas \Beas I_2:=\int_{\mathbb
B^n}|f(z)|\left(|g(z)|+|\nabla g(z)|\right)|\nabla
b(z)|(1-|z|^2)dV(z),\Eeas  and \Beas I_3:=\int_{\mathbb
B^n}|g(z)||\nabla f(z)||\nabla b(z)|(1-|z|^2)dV(z).\Eeas
For the first one, we observe that since $g$ and $b$ are in
all $\mathcal H^p(\mathbb B^n)$, the estimate
$|g(z)b(z)|\le C(1-|z|^2)^{-1/2}$ holds. It follows using
the fact that the measure $(1-|z|^2)^{-1/2}dV(z)$ is a
Carleson measure that
$$I_1\le C\int_{\B^n}| f(z)|(1-|z|^2)^{-1/2}dV(z)\le
C||f||_1.$$ For $I_2$, we use Schwarz inequality to obtain
$$I_2^2\le C\int_{\mathbb
B^n}|f(z)|\left(|g(z)|^2+|\nabla
g(z)|^2\right)|(1-|z|^2)dV(z)\times \int_{\mathbb
B^n}|f(z)||\nabla b(z)|^2(1-|z|^2)dV(z).$$ We conclude by
using the fact that $|\nabla g(z)|^2(1-|z|^2)dV(z)$,
$|\nabla b(z)|^2(1-|z|^2)dV(z)$ and $|
g(z)|^2(1-|z|^2)dV(z)$ are Carleson measures.

The main point is the estimate of $I_3$. We first recall
that, by the weak factorization theorem (see \cite{CRW}), any
$f\in \mathcal H^1(\B^n)$ can be written as
$$f=\sum_jh_jl_j\,\,\,\,\, \textrm{with}\,\,\,\,\sum_j||h_j||_2||l_j||_2\le C||f||_1.$$
Replacing $f$ by this weak factorization, we are led to
estimate a sum of terms as
$$J:=\int_{\B^n}|g(z)||l(z)||\nabla h(z)||\nabla b(z)|(1-|z|^2)dV(z)$$
for  $l$ and $h$ in $\mathcal H^2(\mathbb B^n)$. We recall
that, for $h\in \mathcal H^2(\mathbb B^n)$,
$$\int_{\B^n}|\nabla h(z)|^2(1-|z|^2)dV(z)\le C||h||_2.$$ Using this last inequality, Schwarz
  Inequality and Theorem (referer au chapitre precedent), we obtain
  \Beas
J&\leq &C||h||_{\mathcal H^2}
\left(\int_{\B^n}|g(z)|^2|l(z)|^2|\nabla b(z)|^2(1-|z|^2)dV(z)\right)^{1/2}\\
&\le& C||g||_{BMOA}||l||_2||h||_2.\Eeas This compltes the
proof of the theorem.\ProofEnd

\section{$(\rho,s)$-Carleson measures with $s>1$}
We consider in this section the case of $(\rho,s)$-Carleson
measures when $s>1$. Using  Theorem \ref{C} and the
following equivalence for the norm of elements of the Bloch
space $\mathcal B$:
$$||f||_{\mathcal B}\approx ||f\circ
\varphi_a-f(a)||_{p,\alpha},\,\,\,0<p<\infty\,\,\,
\textrm{and}\,\,\, \alpha>-1$$ (see \cite{KZ}) we can prove
in the same way as Theorem \ref{main1}, the following
theorem.

\btheo\label{main2} Let $0\le p,q<\infty$, $1<s<\infty$.
Let $\mu$ be a positive Borel measure on $\mathbb B^n$.
Then the following conditions are
equivalent.\begin{itemize}\item[i)] There is $C_1
>0$ such that for any $\xi\in \mathbb S^n$ and  $0< \delta <1$,
$$\mu(Q_\delta(\xi))\le C_1\frac{(\sigma(B_\delta(\xi)))^s}{(\log\frac{4}{\delta})^p(\log\log\frac{e^4}{\delta})^q}.$$
\item[ii)] There is $C_2 >0$ such that $$\sup_{a\in \mathbb
B^n}(\log\frac{4}{1-|a|})^p(\log\log\frac{e^4}{1-|a|})^q\int_{\mathbb
B^n}K_a(z)^sd\mu(z)\le C_2 <\infty.$$\item[iii)]There is
$C_3 >0$ such that for any $f\in \mathcal B$,  $$\sup_{a\in
\mathbb
B^n}(\log\log\frac{e^4}{1-|a|})^q\int_{B_n}|f(z)|^{p}K_{a}^s(z)d\mu(z)\le
C_3||f||_{\mathcal B}^{p}.$$\item[iv)]There is $C_4 >0$
such that for any $g\in L\mathcal B$,  $$\sup_{a\in \mathbb
B^n}(\log\frac{4}{1-|a|})^p\int_{\mathbb
B^n}|g(z)|^{q}K_{a}^s(z)d\mu(z)\le C_4||f||_{L\mathcal
B}^{q}.$$\item[v)]There is $C_5 >0$ such that for any $f\in
\mathcal B$ and any $g\in L\mathcal B$,  $$\sup_{a\in
\mathbb B^n}\int_{\mathbb
B^n}|f(z)|^p|g(z)|^{q}K_{a}(z)^sd\mu(z)\le
C_5||f||_{\mathcal B}^{p}||g||_{L\mathcal B}^{q}.$$

\item[vi)] For $0<r<\infty$, there is $C_6
>0$ such that for any $f\in \mathcal B$ and any $g\in L\mathcal B$ and any $h\in A_{ns-(n+1)}^r(\mathbb B^n)$,
$$\int_{\mathbb B^n}|f(z)|^{p}|g(z)|^{q}|h(z)|^rd\mu(z)\le
C_6||f||_{\mathcal B}^{p}||g||_{L\mathcal
B}^{q}||h||_{ns-(n+1),r}^r.$$\end{itemize}\etheo

We now move to applications of Theorem \ref{main2}. We
begin by considering the boundedness of the operator $T_b$
on the logarithmic Bloch space $L\mathcal {B}$. It is not
hard to see that a function $f\in H(\mathbb B^n)$ is in
$L\mathcal B$ if and only if for any $s>1$ the measure
$(1-|z|^2)^{n(s-1)+1}|Rf(z)|^2dV(z)$ is $(\rho,s)$-Carleson
measure with $\rho(t)=(\log(4/t))^2$ (see also Lemma
\ref{logBloch} below) or equivalently that
$$\sup_{a\in
\mathbb B^n}(\log\frac{4}{1-|a|})^2\int_{\mathbb
B^n}K_{a}^{s}(z)(1-|z|^2)^{n(s-1)+1}|Rf(z)|^2dV(z)<\infty.$$
We have the following corollary. \bcor For $b\in H(\mathbb
B^n)$, the operator $T_b$ is bounded on $L\mathcal {B}$ if
and only for any $s>1$, \Be\label{T_b4}I_b<\infty,\Ee where
\Beas I_b &=& \sup_{a\in \mathbb
B^n}(\log\frac{4}{1-|a|})^2(\log\log\frac{e^4}{1-|a|})^2\int_{\mathbb
B^n}|Rb(z)|^2(1-|z|^2)^{n(s-1)+1}K_{a}^{s}(z)dV(z).\Eeas
\ecor\Proof Let $$J_b(f)=\sup_{a\in \mathbb
B^n}(\log\frac{4}{1-|a|})^2\int_{\mathbb
B^n}K_{a}^{s}(z)(1-|z|^2)^{n(s-1)+1}|f(z)|^2|Rb(z)|^2dV(z).$$That
$T_b$ is bounded on $L\mathcal {B}$ is equivalent to saying
there exists a constant $C>0$ such that for any $s>1$ and
any $f\in L\mathcal {B}$,
$$J_b(f)<C||f||_{L\mathcal B}^2$$
which by Theorem \ref{main2} is equivalent to
(\ref{T_b4}).\ProofEnd

Using Theorem \ref{main2} and the fact that any holomorphic
function $f$ belongs to $\mathcal B$ if and only if the
measure $|Rf(z)|^2(1-|z|^2)^{n(s-1)+1}dV(z)$ is a
s-Carleson measure for any $s>1$, we can prove the
following result in the same way.\bcor For $b\in H(\mathbb
B^n)$, the operator $T_b$ is bounded from $L\mathcal {B}$
to $\mathcal {B}$ if and only for $s>1$
\Be\label{T_b5}\sup_{a\in \mathbb
B^n}(\log\log\frac{e^4}{1-|a|})^2\int_{\mathbb
B^n}|Rb(z)|^2(1-|z|^2)^{n(s-1)+1}K_{a}^{s}(z)dV(z)<\infty.\Ee
\ecor

The following well-known result (see for example \cite{S})
follows the same way.\bcor For $b\in H(\mathbb B^n)$, the
operator $T_b$ is bounded on $\mathcal {B}$ if and only for
$s>1$ \Be\label{T_b6}\sup_{a\in \mathbb
B^n}(\log\frac{4}{1-|a|})^2\int_{\mathbb
B^n}|Rb(z)|^2(1-|z|^2)^{n(s-1)+1}K_{a}^{s}(z)dV(z)<\infty.\Ee
\ecor

We also obtain as in the previous section the following
characterization of multipliers of Bloch-type
spaces.\bcor\label{LBlockMultip} An analytic function $f$
on $\mathbb B^n$ belongs to $\mathcal{M}(L\mathcal B)$ if
and only if $f\in \mathcal H^\infty(\mathbb B^n)$ and
satisfies (\ref{T_b4}).\ecor\bcor\label{LBlockBlockMultip}
An analytic
function $f$ on $\mathbb B^n$ belongs to\\
$\mathcal{M}(L\mathcal B,\mathcal B)$ if and only if $f\in
H^\infty(\mathbb B^n)$ and satisfies
(\ref{T_b5}).\ecor\bcor\label{BlockMultip}
An analytic function $f$ on $\mathbb B^n$ belongs to\\
$\mathcal{M}(\mathcal B)$ if and only if $f\in \mathcal
H^\infty(\mathbb B^n)$ and satisfies (\ref{T_b6}).\ecor

\section{Some generalizations}

We give some generalizations and their applications. The
proofs here follow the same steps as in the two previous
sections.

\btheo\label{main3} Let $0\le p_1,p_2,q_1,q_2<\infty$ and
let $\mu$ be a positive Borel meaasure on $\mathbb B^n$.
Then the following conditions are
equivalent.\begin{itemize}\item[i)] There is $C_1
>0$ such that for any $\xi\in \mathbb S^n$ and  $0< \delta <1$, $$\mu(Q_\delta(\xi))\le
C_1\frac{\sigma(B_\delta(\xi))}{(\log\frac{4}{\delta})^{p_1+p_2}(\log\log\frac{e^4}{\delta})^{q_1+q_2}}.$$
\item[ii)]There is $C_2 >0$ such that for any $f\in BMOA$
and any $g\in LMOA$ $$I(f,g) \le
C_2||f||_{BMOA}^{p_1}||g||_{LMOA}^{q_1},$$ \Beas I(f,g) &=&
\sup_{a\in \mathbb
B^n}(\log\frac{4}{1-|a|})^{p_2}(\log\log\frac{e^4}{1-|a|})^{q_2}\int_{\mathbb
B^n}|f_1(z)|^{p_1}|g_1(z)|^{q_1}K_{a}(z)d\mu(z).\Eeas\item[iii)]There
is $C_3
>0$ such that for any $g\in LMOA$
$$I(g)\le C_3||g||_{LMOA}^{q_1},$$
$$I(g)=\sup_{a\in \mathbb
B^n}(\log\frac{4}{1-|a|})^{p_1+p_2}(\log\log\frac{e^4}{1-|a|})^{q_2}\int_{\mathbb
B^n}|g(z)|^{q_1}K_{a}(z)d\mu(z).$$ \item[iv)]There is $C_4
>0$ such that for any $f\in BMOA$
$$I(f)\le C_4||f||_{BMOA}^{p_1},$$
$$I(f)=\sup_{a\in \mathbb
B^n}(\log\frac{4}{1-|a|})^{p_2}(\log\log\frac{e^4}{1-|a|})^{q_1+q_2}\int_{\mathbb
B^n}|f(z)|^{p_1}K_{a}(z)d\mu(z).$$
\end{itemize}\etheo\btheo\label{main4} Let $0\le p_1,p_2,q_1,q_2<\infty$, let $1<s<\infty$ and
$\mu$ be a positive Borel measure on $\mathbb B^n$. Then
the following conditions are
equivalent.\begin{itemize}\item[i)] There is $C_1
>0$ such that for any $\xi\in \mathbb S^n$ and  $0< \delta <1$, $$\mu(Q_\delta(\xi))\le
C_1\frac{(\sigma(B_\delta(\xi)))^s}{(\log\frac{4}{\delta})^{p_1+p_2}(\log\log\frac{e^4}{\delta})^{q_1+q_2}}.$$
\item[ii)]There is $C_2 >0$ such that for any $f\in
\mathcal B$ and any $g\in L\mathcal B$,
$$J(f,g)\le
C_2||f||_{\mathcal B}^{p_1}||g||_{L\mathcal B}^{q_1},$$
\Beas J(f,g) &=& \sup_{a\in \mathbb
B^n}(\log\frac{4}{1-|a|})^{p_2}(\log\log\frac{e^4}{1-|a|})^{q_2}\int_{\mathbb
B^n}|f_1(z)|^{p_1}|g_1(z)|^{q_1}K_{a}^s(z)d\mu(z).\Eeas\item[iii)]There
is $C_3
>0$ such that for any $g\in L\mathcal B$,
$$\sup_{a\in \mathbb B^n}(\log\frac{4}{1-|a|})^{p_1+p_2}(\log\log\frac{e^4}{1-|a|})^{q_2}\int_{\mathbb
B^n}|g(z)|^{q_1}K_{a}^s(z)d\mu(z)\le C_3||g||_{L\mathcal
B}^{q_1}.$$ \item[iv)]There is $C_4
>0$ such that for any $f\in \mathcal B$
$$\sup_{a\in \mathbb B^n}(\log\frac{4}{1-|a|})^{p_2}(\log\log\frac{e^4}{1-|a|})^{q_1+q_2}\int_{\mathbb
B^n}|f(z)|^{p_1}K_{a}^s(z)d\mu(z)\le C_4||f||_{\mathcal
B}^{p_1}.$$
\end{itemize}\etheo

Let $0\le p,q<\infty$. An analytic function $f$ belongs to
$BMOA_{\rho_{p,q}}$ with
$\rho_{p,q}(t)=(\log(4/t))^p(\log\log(e^4/t))^q$ if
 $f\in \mathcal H^1(\mathbb B^n) $ and there exists a constant $C>0$ so that
$$\sup_{B=B_\delta(\xi)\atop \delta\in]0,1[,\xi\in \mathbb S^n} \frac {(\log(4/\delta))^p(\log\log(e^4/\delta))^q}{\sigma(B)}\int _B|f-f_{B}|d\sigma\le  C.$$
 By \cite{WS}, a function $f$
belongs to $BMOA_{\rho_{p,q}}$ if and only if
$d\mu(z)=(1-|z|^2)|\nabla f(z)|^2dV(z)$ is a
$\rho_{p,q}^2$- Carleson measure. The following corollaries
can be proved as in the previous sections.\bcor Let $0\le
p,q<\infty$. Given an analytic function $b$, the operator
$T_b$ is bounded from $LMOA$ to $BMOA_{\rho_{p,q}}$ if and
only if \Be\label{T_b7}\sup_{a\in \mathbb
B^n}(\log\frac{4}{1-|a|})^{2p}(\log\log\frac{e^4}{1-|a|})^{2q+2}\int_{\mathbb
B^n}|Rb(z)|^2(1-|z|^2)K_{a}(z)dV(z)<\infty.\Ee\ecor

\bcor Let $0\le p,q<\infty$. Given an analytic function
$b$, the operator $T_b$ is bounded from $BMOA$ to
$BMOA_{\rho_{p,q}}$ if and only if
\Be\label{T_b8}\sup_{a\in \mathbb
B^n}(\log\frac{4}{1-|a|})^{2p+2}(\log\log\frac{e^4}{1-|a|})^{2q}\int_{\mathbb
B^n}|Rb(z)|^2(1-|z|^2)K_{a}(z)dV(z)<\infty.\Ee\ecor

In particular, we have the following.\bcor  Given an
analytic function $b$, the operator $T_b$ is bounded from
$BMOA$ to $LMOA$ if and only if \Beas \sup_{a\in \mathbb
B^n}(\log\frac{4}{1-|a|})^{4}\int_{\mathbb
B^n}|Rb(z)|^2(1-|z|^2)K_{a}(z)dV(z)<\infty.\Eeas\ecor
Let us now move to applications of Theorem \ref{main4}. We
first introduce the following generalized
$\alpha$-logarithmic-type Bloch spaces.
\begin{definition} For $0\le p,q<\infty$ and $\alpha>0$. Let denote $\mathcal B_{\alpha}^{p,q}$
the space of holomorphic functions $f$ such that
$$\sup_{z\in \mathbb
B^n}(1-|z|^2)^\alpha|Rf(z)|(\log\frac{4}{1-|z|^2})^p(\log\log\frac{e^4}{1-|z|^2})^q<\infty.$$
\end{definition} These can be seen as special case of the so called
$\mu$-Bloch spaces (see for example \cite{HZJ}) and one
have that $\mathcal B_\alpha^{p,q}$ are Banach spaces with
the norm \Beas ||f||_{\mathcal B_\alpha^{p,q}} &=&
|f(0)|+\sup_{z\in \mathbb
B^n}(1-|z|^2)^{\alpha}|Rf(z)|(\log\frac{4}{1-|z|^2})^p(\log\log\frac{e^4}{1-|z|^2})^q<\infty.\Eeas
The usual Bloch space $\mathcal B$ then corresponds to the
case $\alpha=1$ and $p=q=0$ while $\mathcal
B_\alpha^{0,0}=\B_\alpha$ are the so called $\alpha$- Bloch
spaces (see \cite{KZ} ) and $\mathcal B_1^{1,1}=L\mathcal
B$ is the usual logarithmic Bloch space.

Let $\varphi_z$ be the involutive automorphism of $\mathbb
B^n$ that interchanges $0$ and $z$. The Bergman metric of
$\mathbb B^n$ is given by
$$\beta(z,w):=\frac{1}{2}\log\frac{1+|\varphi_z(w)|}{1-|\varphi_z(w)|},$$
for all $z,w\in \mathbb B^n$. For any $R>0$ and any $a\in
\mathbb B^n$, we write
$$D(a,R)=\{z\in \mathbb B^n: \beta(z,a)<R\}$$ for the Bergman ball centered at $a$ with
radius $R$.
We have the following characterization of elements of
$\mathcal B_{\alpha}^{p,q}$.
 \blem\label{logBloch} Let $0\le p,q<\infty$ and $\alpha>0$. A
function $f\in H(\mathbb B^n)$ is in $\mathcal
B_{\alpha}^{p,q}$ if and only if
$(1-|z|^2)^{n(s-1)+2\alpha-1}|Rf(z)|^2dV(z)$ is a
$(\rho_{p,q},s)$-Carleson measure for any $s>1$, where
$\rho_{p,q}(t)=(\log\frac{4}{t})^{2p}(\log\log\frac{e^4}{t})^{2q}$.\elem
\Proof  We first suppose that $f$ belongs to $\mathcal
B_{\alpha}^{p,q}$ and show that there exists a constant
$C>0$ such that for any $\xi\in \mathbb S^n$, $0<\delta<1$
and any $s>1$, the following inequality holds
\Beas\label{logequiv}I_f(\delta)\le
C\sigma(B_\delta(\xi))^{s},\Eeas where
$$I_f(\delta)=(\log\frac{4}{\delta})^{2p}(\log\log\frac{e^4}{\delta})^{2q}\int_{Q_\delta(\xi)}|Rf(z)|^2(1-|z|^2)^{n(s-1)+2\alpha-1}dV(z).$$
Let
$h(x)=(\log\frac{4}{x})^{2p}(\log\log\frac{e^4}{x})^{2q}$,
$h$ is decreasing on $(0,1)$ and moreover, for any $z\in
Q_\delta(\xi)$, $1-|z|^2<|1-<\xi,z>|<\delta$. It follows
using the definition of $\mathcal B_{\alpha}^{p,q}$ that
there exists a constant $C>0$ so that for all $f\in
\mathcal {B}_\alpha^{p,q}$, \Beas I_f(\delta) &\le&
C\int_{Q_\delta(\xi)}\frac{h(\delta)}{h(1-|z|^2)}(1-|z|^2)^{n(s-1)-1}dV(z)\\
&\le& C\int_{Q_\delta(\xi)}(1-|z|^2)^{n(s-1)-1}dV(z)\\
&\le& C\sigma(B_\delta(\xi))^{s}.\Eeas This shows the
necessary part.

Conversely, let us suppose that the analytic function $f$
is such that there exists $C>0$ so that for any $\xi\in
\mathbb S^n$, $0<\delta<1$ and any $s>1$,$$I_f(\delta)\le
C(\sigma(B_\delta(\xi)))^{s}$$ and show that in this case
$f$ belongs to $\mathcal B_\alpha^{p,q}$. We recall that
for $a\in \mathbb B^n$, letting $\delta=1-|a|$, for $R>0$,
there exists $\lambda\in (0,1)$ such that $D(a,R)\subset
Q_\delta(\xi)$ with $a=(1-\lambda\delta)\xi$ (see Lemma
5.23 of \cite{KZ}). Now, using the mean value property, we
obtain that for any $a\in \mathbb B^n$, $$|Rf(a)|^2\le
\frac{C}{(1-|a|^2)^{ns+2\alpha}}\int_{D(a,R)}|Rf(z)|^2(1-|z|^2)^{n(s-1)+2\alpha-1}dV(z).$$
It follows from the above inclusion and the hypotheses on
the measure\\ $|Rf(z)|^2(1-|z|^2)^{n(s-1)+2\alpha-1}dV(z)$
that $$
(1-|a|^2)^{2\alpha}|Rf(a)|^2(\log\frac{4}{1-|a|^2})^{2p}(\log\log\frac{e^4}{1-|a|^2})^{2q}
\le \frac{C}{\delta^{ns}}I_f(\delta)\le C<\infty.$$  The
proof is complete.\ProofEnd

The following two corollaries can be proved exactly as
before.\bcor\label{Blochtologbloch} Let $0\le p,q<\infty$,
$\alpha>0$ and $b\in H(\mathbb B^n)$. Then the following
conditions are equivalent.\begin{itemize} \item[(a)]$T_b$
is bounded from $L\mathcal {B}$ to $\mathcal
B_\alpha^{p,q}$.\item[(b)] For any $s>1$ \Be\label{T_b9}
I_{b,p,q}<\infty,\Ee \Beas I_{b,p,q} &=& \sup_{a\in
B_n}(\log\frac{4}{1-|a|})^{2p}(\log\log\frac{e^4}{1-|a|})^{2q+2}\int_{\mathbb
B^n}|Rb(z)|^2(1-|z|^2)^{n(s-1)+2\alpha-1}K_{a}^s(z)dV(z).\Eeas\end{itemize}
\ecor\bcor Let $0\le p,q<\infty$, $\alpha>0$ and $b\in
H(\mathbb B^n)$. Then the following conditions are
equivalent.\begin{itemize} \item[(a)]$T_b$ is bounded from
$\mathcal {B}$ to $\mathcal B_\alpha^{p,q}$.\item[(b)] For
any $s>1$ \Be\label{T_b10}J_{b,p,q}<\infty,\Ee \Beas
J_{b,p,q} &=& \sup_{a\in
B_n}(\log\frac{4}{1-|a|})^{2p+2}(\log\log\frac{e^4}{1-|a|})^{2q}\int_{\mathbb
B^n}|Rb(z)|^2(1-|z|^2)^{n(s-1)+2\alpha-1}K_{a}^s(z)dV(z).\Eeas\end{itemize}\ecor
In particular, we have the following.\bcor Given $b\in
H(\mathbb B^n)$, the operator $T_b$ is bounded from
$\mathcal {B}$ to $L\mathcal {B}$ if and only if for any
$s>1$, \Beas\sup_{a\in
B_n}(\log\frac{4}{1-|a|})^{4}\int_{\mathbb
B^n}|Rb(z)|^2(1-|z|^2)^{n(s-1)+2\alpha-1}K_{a}^s(z)dV(z)<\infty.\Eeas\ecor

\bibliographystyle{plain}

\end{document}